\documentclass{article} 
\usepackage{graphicx,amsmath, amssymb, verbatim}
\usepackage{framed}
\usepackage{fancybox}
\usepackage{amsthm}
\usepackage[all]{xy}
\usepackage{bm}
\usepackage{ascmac}
\usepackage{physics}
\usepackage{color}
\usepackage{comment}
\allowdisplaybreaks[2]
\newtheorem{proposition}{Proposition}
\newtheorem{theorem}{Theorem}

\newtheorem*{definition*}{Definition}

\newcommand {\C } {\mathbb{C}} 

\newcommand {\R } {\mathbb{R}}

\newcommand{\ba}{\begin{eqnarray}}
\newcommand{\ea}{\end{eqnarray}}

\begin{document} 

\title{\bf{ Schubert Calculus via Fermionic Variables}} 
\author{\large  Ken Kuwata \\
\\ \it Department of General Education \\
\it National Institute of Technology, Kagawa College \\
\it  Chokushi, Takamatsu, 761-8058, Japan\\
\\
\it e-mail address: kuwata-k@t.kagawa-nct.ac.jp  }

\maketitle 

\begin {abstract} 
Imanishi, Jinzenji and Kuwata provided a recipe for computing Euler number of Grassmann manifold $G(k,N)$ using physical model and its path-integral [S.Imanishi, M.Jinzenji and K.Kuwata, Journal of Geometry and Physics, Volume 180, October 2022, 104623]. They demonstrated that the cohomology ring of $G(k,N)$ is represented by fermionic variables. In this study, using only fermionic variables, we computed an integral of the Chern classes of the dual bundle of the tautological bundle on $G(k,N)$. In other words, the intersection number of  the Schubert cycles is obtained using the fermion integral. 
\end {abstract}

\section{Introduction }
\subsection{Background}
In this study, we aim to compute the intersection numbers of Schubert cycles. We used fermionic variables and their integrals in \cite{IJK}. In this section, we explain the background of the study. 
The complex Grassmann manifold $G(k,N)$ is the space parameterizing all  $k$ -dimensional linear subspaces of $N$-dimensional complex vector space $\C^N$. Because the elements of its cohomology ring are represented by the Poincar\'{e} dual of some Schubert cycles of $G(k,N)$, their integral provides the intersection number of  Schubert cycles. This research is called Schurbert calculus, and has been studied in combinatorics, representation theory, and other fields \cite{IN}. The integral of these cohomology classes can be computed using localization theory or the Landau-Ginzburg formulation. In the localization theory, a fixed-point theorem for a compact manifold with torus action is used. In particular, the formula for the intersection number is provided using the localization theory \cite{AW,DTH, MZ}. 
However, the Landau-Ginzburg formulation \cite{EW3, NC} uses a potential function provided by the total Chern class of the tautological bundle of $G(k,N)$  and residue. However, we do not use these theories. We employed the theory of \cite{IJK}.  Imanishi et al. constructed a physical toy model for computing the Euler number of  $G(k,N)$. The model was constructed using two types of variables. One is a commutative variable called a bosonic variable, while the other is an anticommutative variable, called a fermionic variable. In \cite{IJK}, it was found that the cohomology ring of $G(k,N)$ can be represented by fermionic variables, and that the Euler number is provided by their integral. Therefore, the intersection number of Schubert cycles can be obtained using fermion integrals. Generally, it is difficult to perform this calculation. However, in some cases, the number of intersections can be calculated using this method. In this study,  we demonstrated the use of the method of \cite{IJK}.

\subsection{Organization of the paper}
This paper is divided into two sections.

In Section 1, we describe our background and theorem. In addition to the background described above, we introduce the relationship between Chern classes  and Schubert cycles, our  theorem in this paper, and the theory in  \cite{IJK}.
First, we remark on Chern classes and Schubert cycles. Next, we introduce the theorem. Finally,  we introduce the relation between the Chern classes and fermionic variables in \cite{IJK}.

In Section 2, we provide the proof of our theorem. We computed the fermion integral to prove the theorem.  We also summarize the important results of the fermion integrals.

\subsection{Chern classes and Schubert cycles}
In this section, we explain the relation between the Chern classes and Schubert cycles and our theorem. In this study, we employed the notation in \cite{IJK}. First, we introduce the cohomology ring of $G(k,N)$ \cite{BT, IJK}, and then remark on a tautological bundle $S$ and a universal quotient bundle $Q$.  The fiber of $S$ at $\Lambda\in G(k,N)$ is the complex $k$-dimensional subspace $\Lambda\subset \mathbb{C}^N$ itself ($\mathrm{rk} (S)=k$). Subsequently, a universal quotient bundle $Q$ ( $\mathrm{rk}(Q)=N-k$) is defined by the following exact sequence:
\begin{align}
0\to S \to \mathbb{C}^N  \to Q \to 0, \label{exact 1}
\end{align}
where $\mathbb{C}^N$ denotes the trivial bundle $G(k,N)\times \mathbb{C}^N$. We write  $c_i(E)$ as the $i$-th Chern class of the vector bundle $E$. Let $E^*$ be the dual bundle of $E$. Subsequently, the cohomology ring $H^{*}(G(k,N))$  of $G(k,N)$ is
\begin{align}
H^{*}(G(k,N))=\frac{\R[c_1(S^*),\cdots,c_k(S^*),c_1(Q^*),\cdots,c_{N-k}(Q^*)]}{(c(S^*)c(Q^*)=1)}. 
\label{crgd}
\end{align} 
Consider the problem of representing  $H^{*}(G(k,N))$ using $c_j(S^*) (j=1,2,\cdots,k)$. This was presented in \cite{IJK}. If we decompose $S^*$ formally by line bundle $L_i$ ($i=1,2,\cdots,k$):
\ba
S^{*}=\mathop{\oplus}_{i=1}^{k}L_{i},
\ea
$c(S^*)$ and $c_i(S^*)$ are expressed as follows: 
\begin{align}
c(S^*)=\prod_{i=1}^k(1+tx_i)
=1+\sum_{j=1}^kt^jc_j(S^*),\;\;(x_i:= c_1(L_i)). \label{decS*}
\end{align} 
Hence, $c_j(S^*)$ is written as the degree $j$ elementary symmetric polynomial of $x_1,\cdots,x_k$. Then, the relation $c(S^*)c(Q^*)=1$ can be rewritten as 
\begin{align}
&c(Q^*)=\frac{1}{c(S^*)}
=\frac{1}{1+\sum_{j=1}^kt^jc_j(S^*)}=\sum_{i=0}^{\infty}a_{i}t^{i}.
\label{grrel}
\end{align}
We can rewrite $a_{i}$ in (\ref{grrel}) as 
\ba
c_{i}(Q^{*})=a_{i}\;\;(i=1,2,\cdots,N-k),\;\;a_{i}=0\;\;(i>N-k).
\ea
Moreover, $a_{i}$ is the degree $i$ homogeneous polynomial of $c_j(S^*)$'s $(j=1,2,\cdots,k)$.  
Thus, $c_{i}(Q^{*})$ in (\ref{crgd}) can be rewritten as $c_i(S^*)$.  Consequently, we obtain another representation of $H^{*}(G(k,N))$:
\begin{align}
H^{*}(G(k,N))=\frac{\R[c_1(S^*),\cdots,c_k(S^*)]}{(a_{i}=0\;\;(i>N-k))}. 
\label{crg2}
\end{align} 
 Second, we introduce the Schubert cycle and explain the relationship between   Chern classes and Schubert cycles. For a more detailed discussion, please refer to \cite{GH}. For any flag $V: 0\subset V_1\subset V_2\subset \cdots \subset V_N=\mathbb{C}^N$, Schubert manifold $\sigma_a(V)$ is defined as follows:
\begin{align}
\sigma_a(V):=\{ \Lambda \in G(k,N)|\mathrm{dim}(\Lambda \cap V_{N-k+i-a_i}) \geq i (1\leq i \leq k)\},
\end{align}
where $a=(a_1,\cdots,a_k)$ denotes a sequence of natural numbers that satisfies $0\leq a_k\leq a_{k-1} \leq \cdots \leq a_1 \leq N-k$. $\sigma_a(V)$ is a subvariety of $G(k,N)$ of dimension $\sum_{l=1}^k a_l$. The homology class of $\sigma_a(V)$ is independent of the chosen flag. Therefore, let  $\sigma_a(V)$ as the homology class be denoted by $\sigma_a$.
Let $\sigma_a^\vee$ be the Poincar\'{e} dual of the cycle $\sigma_a$. For  simplicity of notation, we omit $0$ from $a$. For example, $\sigma_{a_1,a_2,\cdots, a_n}$ denotes  $\sigma_{(a_1,a_2,\cdots, a_n,0,\cdots,0)}$. The relationship between $i$-th Chern class of a vector bundle $E$ and that  of its dual bundle $E^*$
is provided by $c_i(E^*)=(-1)^ic_i(E)$. From this formula and the Gauss-Bonnet theorem, we obtain: 
\begin{align}
c_i(S^*)=(-1)^ic_i(S)=\sigma^\vee_{\footnotesize{\underbrace{1,\cdots,1}_{\normalsize{i}}}}=:\sigma_{1^{(i)}}^\vee. 
\end{align}
Finally, we introduce. 
\begin{theorem}\label{theo1}
\begin{align}
\int_{G(k,N)}(\sigma^\vee_{1^{(1)}})^{kN-k^2}&=(kN-k^2)!\frac{ \prod_{j=0}^{k-1} j!}{ \prod_{j=N-k}^{N-1} j!} . \label{theo1-1}\\
\int_{G(k,N)} (\sigma^\vee_{1^{(1)}})^{kN-k^2-2}(\sigma^\vee_{1^{(2)}}) 
&= \frac{(kN-k^2-2)!(N-k)(N-k+1)k(k-1)}{2}\frac{ \prod_{j=0}^{k-1} j!}{ \prod_{j=N-k}^{N-1} j!}.  \label{theo1-2}\\
\int_{G(k,N)} (\sigma^\vee_{1^{(1)}})^{kN-k^2-4}(\sigma^\vee_{1^{(2)}})^2 
&= \frac{(kN-k^2-4)!(N-k)(N-k+1)k(k-1)}{4} \frac{\prod_{j=0}^{k-1} j!}{\prod_{j=N-k}^{N-1} j!}\notag \\
&\times \Bigl[k(k-1)(N-k)(N-k-1) +2(k-2)(k-3)(N-k)\notag \\
&+4(k-2)(N-k-1) \Bigr] \label{theo1-3}.
\end{align}
Here, we assume that $N$ and $k$ in (\ref{theo1-2}) and (\ref{theo1-3}) satisfy $kN-k^2-2 \geq 0$ and $kN-k^2-4 \geq 0$, respectively.
\end{theorem} 
Note that these are the intersection numbers of $\sigma_{1^{(1)}}$ and $\sigma_{1^{(2)}}$. 
However, the results of (\ref{theo1-1}) are already well known \cite{NC,WF}. 
 When $k=2$, the intersection numbers of $\sigma^*_{1^{(1)}}$ and $\sigma^*_{1^{(2)}}$ in $G(2,N)$ are known \cite{NC}. 
\begin{align}
\int_{G(2,N)} (\sigma^\vee_{1^{(1)}})^{2N-4-2l}(\sigma^\vee_{1^{(2)}})^l 
&= \frac{(2(N-2-l))! }{(N-2-l)!(N-1-l)!}.
\end{align}
(It is also derived by S. Imanishi's Masters thesis  using fermionic variables \cite{Shoichiro}.)
In particular, we obtain the following results from (\ref{theo1-1}), (\ref{theo1-2}), and (\ref{theo1-3}):
 \begin{align}
 \int_{G(2,N)} (\sigma^\vee_{1^{(1)}})^{2N-4}
&= \frac{(2N-4)! }{(N-2)!(N-1)!}, \\
\int_{G(2,N)} (\sigma^\vee_{1^{(1)}})^{2N-6}(\sigma^\vee_{1^{(2)}}) 
&=  \frac{(2N-6)!}{ (N-3)!(N-2)!}, \\
\int_{G(2,N)} (\sigma^\vee_{1^{(1)}})^{2N-8}(\sigma^\vee_{1^{(2)}})^2 
&= \frac{(2N-8)! }{(N-4)!(N-3)!}.
\end{align}

\subsection{Fermionic variables and Cohomology ring of $G(k,N)$ (Review of \cite{IJK})}
We summarize the representation of the cohomology ring of $G(k,N)$ using fermionic variables \cite{IJK}. 
We introduce the fermionic variables $\psi^j_{s}, \;\psi^{\bar j}_{\bar s}$
$(s=1,\cdots,N-k,\;j=1,\cdots k)$ and $(k\times k)$ matrix
\begin{align}
\Phi&:=\sum_{s=1}^{N-k} \left( \begin{array}{ccc}
\psi_s^{1} \psi_{\bar s}^{\bar1} &\ldots&\psi_s^{1} \psi_{\bar s}^{\bar k}\\
\vdots& \ddots&\vdots\\
\psi_s^{k} \psi_{\bar s}^{\bar1}&\ldots&\psi_s^{k} \psi^{\bar k}_{\bar1}\\
\end{array}
\right).
\end{align}
The fermionic variables $\psi^j_{s}, \;\psi^{\bar j}_{\bar s}$ satisfy the following conditions.
\begin{align}
\psi^j_{s}\psi^j_{s}= \psi_{\bar s}^{\bar j}\psi_{\bar s}^{\bar j}=0,~~
\psi^j_{s}\psi^i_{l}=-\psi^i_{l}\psi^j_{s},~\psi_{\bar s}^{\bar j}\psi_{\bar l}^{\bar i}=-\psi_{\bar l}^{\bar i}\psi_{\bar s}^{\bar j},~
\psi^j_{s}\psi_{\bar l}^{\bar i}=-\psi_{\bar l}^{\bar i}\psi^j_{s}
\end{align}
$(s,l=1,2,\cdots,N-k,~i,j=1,2,\cdots,k)$.
The fermion integral is defined as follows:
\begin{align}
\int D\psi \prod_{s=1}^{N-k} \psi_s^{1} \psi_{\bar s}^{\bar1} \cdots \psi_s^{k} \psi_{\bar s}^{\bar k}=1,
\end{align} 
where $D\psi:=\prod_{s=1}^{N-k} d\psi_s^{1} d\psi_{\bar s}^{\bar1} \cdots d\psi_s^{k} d\psi_{\bar s}^{\bar k}$.  We define $\tau_{j}$ $(j=1,2\cdots,k)$ as
\begin{align}
1+\tau_1t+\cdots+\tau_k t^k:=\det(I_k+t \Phi)=\prod_{j=1}^{k}(1+\lambda_{j}t). \label{Deftau}
\end{align} 
Here, $\lambda_j (j=1,\cdots,k)$ are eigenvalues of $\Phi$. 
Specifically, $\tau_{j}$ is the degree $j$ elementary symmetric polynomial of $\lambda_{1},\cdots,\lambda_{k}$.
Note that $\tau_{k}$ is identified with $\det(\Phi)$ and $\tau_{1}$ is identified with $\tr(\Phi)$. 
In \cite{IJK}, the following theorems were proved:
\begin{theorem}\label{ThV}
\rm{\cite{IJK}}
\begin{align}
\frac{ \prod_{j=0}^{k-1} j!}{ \prod_{j=N-k}^{N-1} j!}\int D\psi \left(\det(\Phi)\right)^{N-k} =1.
\end{align}
\end{theorem}
\begin{theorem}\label{theoHG}\rm{\cite{IJK}}
\begin{align}
H^*(G(k,N))=\frac{\R[c_1 (S^*),\cdots,c_k (S^*)]}{(a_{i}=0\;\;(i>N-k))} \simeq \R [\tau_1,\cdots,\tau_k].
\end{align}
\end{theorem}
Theorem \ref{theoHG} is provided by ring homomorphism $f:\R[c_1 (S^*),\cdots,c_k (S^*)]\rightarrow \R[\tau_{1},\cdots,\tau_{k}]$, which is defined  as 
\ba
f(c_j(S^*))=\tau_{j}\;\;(j=1,2,\cdots,k).
\ea  
From the isomorphism  $H^{*}(G(k,N))\xrightarrow{\sim} \R[\tau_{1},\cdots,\tau_{k}]$, $x_{j}$ is identified as $\lambda_{j}$.
Theorem \ref{ThV} corresponds to the normalization condition of the integration
on $G(k,N)$  given by
\ba 
\int_{G(k,N)}(\sigma^\vee_{1^{(k)}})^{N-k}=1.
\ea
Therefore, we obtain the following formula:
\ba
\int_{G(k,N)}g(x_{1},\cdots,x_{k})=\frac{ \prod_{j=0}^{k-1} j!}{ \prod_{j=N-k}^{N-1} j!}\int D\psi g(\lambda_{1},\cdots,\lambda_{k}),
\label{intrel}
\ea
where $g(x_{1},\cdots,x_{k})$ is a symmetric polynomial of $x_{1},\cdots,x_{k}$ that represents an element of $H^{*}(G(k,N))$.

\section{Proof of our theorem}
\subsection{Proof of Theorem \ref{theo1}}
\begin{proof}(Theorem \ref{theo1})

First, we prove (\ref{theo1-1}).  From (\ref{decS*}), (\ref{Deftau}), and (\ref{intrel}), we have
\begin{align}
&\int_{G(k,N)} \Bigl(\sigma_{1^{(1)}}^\vee \Bigr)^{kN-k^2}
=\frac{ \prod_{j=0}^{k-1} j!}{ \prod_{j=N-k}^{N-1} j!}\int D\psi \left(\mathrm{tr}\left(\Phi \right) \right)^{kN-k^2} 
=\frac{ \prod_{j=0}^{k-1} j!}{ \prod_{j=N-k}^{N-1} j!}\int D\psi \left(\sum_{s=1}^{N-k} \sum_{j=1}^k\psi_s^j \psi_{\bar s}^{\bar j} \right)^{kN-k^2}.
\end{align}
From the multinomial theorem and the conditions of the fermionic variables $\psi_s^j \psi_s^j=\psi_{\bar s}^{\bar j}\psi_{\bar s}^{\bar j}=0$, we obtain:
\begin{align}
\int_{G(k,N)} \Bigl(\sigma_{1^{(1)}}^\vee \Bigr)^{kN-k^2}
&=(Nk-k^2)! \frac{ \prod_{j=0}^{k-1} j!}{ \prod_{j=N-k}^{N-1} j!}\int D\psi \prod_{s=1}^{N-k} \prod_{j=1}^k\psi_s^j \psi_{\bar s}^{\bar j} =(kN-k^2)! \frac{ \prod_{j=0}^{k-1} j!}{ \prod_{j=N-k}^{N-1} j!}.
\end{align}
Second, we show that (\ref{theo1-2}) and (\ref{theo1-3}). 
 In the same way as in (\ref{theo1-1}),
 \begin{align}
 \int_{G(k,N)} (\sigma^\vee_{1^{(1)}})^{kN-k^2-2l}(\sigma^\vee_{1^{(2)}}) ^l
&= \frac{ \prod_{j=0}^{k-1} j!}{ \prod_{j=N-k}^{N-1} j!}  
\int D\psi (\tau_1)^{kN-k^2-2l}(\tau_2)^l ~(l=1,2).
 \end{align}
 As $\tau_2=\frac{1}{2}\{(\mathrm{tr}\left( \Phi \right))^2-\mathrm{tr}\left( \Phi^2 \right)\}$,
\begin{align}
&\int D\psi (\tau_1)^{kN-k^2-2l}(\tau_2)^l
=\frac{1}{2^l}\int D\psi \left\{\mathrm{tr}\left( \Phi \right)\right\}^{kN-k^2-2l} \left\{ (\mathrm{tr}\left( \Phi \right))^2-\mathrm{tr}\left( \Phi^2 \right)\right\}^l \\
&=\frac{1}{2^l} \sum_{m=0}^l \binom{l}{m}(-1)^m \int D\psi \left( \mathrm{tr}\left( \Phi \right)\right)^{kN-k^2-2m} \left(\mathrm{tr}\left( \Phi^2 \right)\right)^m.
\end{align}
Let us define
\begin{align}
P_m:=\int D\psi \left( \mathrm{tr}\left( \Phi \right)\right)^{kN-k^2-2m} \left(\mathrm{tr}\left( \Phi^2 \right)\right)^m~(m=0,1,2).
\end{align}
As can be observed from the calculation in (\ref{theo1-1}), $P_0=(kN-k^2)!$.
We can obtain the following result for $P_1$ and $P_2$.
\begin{proposition}\label{prop1}
\begin{align}
P_1&=(kN-k^2-2)!k(N-k)(N-2k)\label{prop1-1}. \\
P_2&=(kN-k^2-4)!k(N-k)\Bigl [k(N-k)^3-2(N-k)^2(k^2+2) +(N-k)(k^3+10k)-4k^2-2\Bigr]. \label{prop1-2}
\end{align}
\end{proposition}
We will prove these results later in this paper. From Proposition \ref{prop1}, we have  
\begin{align}
\int D\psi (\tau_1)^{kN-k^2-2}(\tau_2)&=\frac{1}{2}(P_0-P_1)
=\frac{1}{2}(kN-k^2-2)!k(N-k)\{(kN-k^2-1)-(N-2k)\} \notag \\
&=\frac{1}{2}(kN-k^2-2)!(N-k)(N-k+1)k(k-1).
\end{align}
We obtain (\ref{theo1-2}). Similarly, we obtain (\ref{theo1-3}) from $\int D\psi (\tau_1)^{kN-k^2-4}(\tau_2)^2=\frac{1}{4}(P_0-2P_1+P_2)$. We have proved Theorem \ref{theo1}.
\end{proof}

\subsection{Proof of Proposition \ref{prop1}}
\begin{proof}
(Proposition \ref{prop1}).

We compute $P_1$.
Let $\omega^{i j}$ be $\sum_{s=1}^{N-k} \psi_s^i \psi_{\bar s}^{\bar j}$.
By definition,
\begin{align}
&\mathrm{tr}\left( \Phi^2 \right)=\sum_{i,j=1}^k\omega^{ij}\omega^{ji}=\sum_{i=1}^k (\omega^{ii})^2+\sum_{i\neq j} \omega^{ij}\omega^{ji}, \\
&P_1=\int D\psi \left( \mathrm{tr}\left( \Phi \right)\right)^{kN-k^2-2} \left(\mathrm{tr}\left( \Phi^2 \right)\right) \\
&=\sum_{i=1}^k \int D\psi \left( \sum_{n=1}^k \omega^{nn} \right)^{kN-k^2-2} (\omega^{ii})^2+\sum_{i\neq j} \int D\psi \left( \sum_{n=1}^k \omega^{nn} \right)^{kN-k^2-2} \omega^{ij}\omega^{ji} \\
&=\sum_{i=1}^k \sum_{p_n}\frac{(kN-k^2-2)!}{\prod_{n=1}^k p_n!} \int D\psi \left( \prod_{n=1}^k (\omega^{nn})^{p_n} \right) (\omega^{ii})^2 \notag \\
&+\sum_{i\neq j} \sum_{p_n}\frac{(kN-k^2-2)!}{\prod_{n=1}^k p_n!} \int D\psi \left( \prod_{n=1}^k (\omega^{nn})^{p_n} \right) \omega^{ij}\omega^{ji}.
\end{align}
Here, $\sum_{p_n}$ indicates that the sum includes all combinations from $0$ to $kN-k^2-2$ indices $p_1$ through $p_k$, such that the sum of all $p_n (n=1,\cdots,k)$ is $kN-k^2-2$. In the first term, because each $\omega^{ii}$ $(i=1,\cdots,k)$ must be $N-k$ for  the fermion integral to be non-zero, $p_n=N-k  (n\neq i)$ and $p_i=N-k-2$. In the second term, $p_n=N-k (n\neq i,j)$ and $p_i=p_j=N-k-1$.
\begin{align}
&P_1=\sum_{i=1}^k \frac{(kN-k^2-2)!}{((N-k)!)^{k-1}(N-k-2)!} \int D\psi \left( \prod_{n=1}^k (\omega^{nn})^{N-k} \right)  \notag \\
&+\sum_{i\neq j} \frac{(kN-k^2-2)!}{((N-k)!)^{k-2}((N-k-1)!)^2} \int D\psi \left( \prod_{n\neq i,j} (\omega^{nn})^{N-k} \right) \left( \omega^{ii}\omega^{jj}\right)^{N-k-1} \omega^{ij}\omega^{ji}.
\end{align}
From $\omega^{ii}=\sum_{s=1}^{N-k} \psi^i_s \psi^{\bar i}_{\bar s}$, the multinomial theorem and conditions of the fermionic variables $\psi_s^j \psi_s^j=\psi_{\bar s}^{\bar j}\psi_{\bar s}^{\bar j}=0$.
\begin{align}
&P_1=\sum_{i=1}^k \frac{(kN-k^2-2)!}{(N-k-2)!}(N-k)!  \notag \\
&+\sum_{i\neq j} \frac{(kN-k^2-2)!}{((N-k-1)!)^2} \int D\psi \left( \prod_{n\neq i,j} \prod_{l=1}^{N-k} \psi^n_l \psi^{\bar n}_{\bar l} \right) \left( \omega^{ii}\omega^{jj}\right)^{N-k-1} \left( \sum_{s,t=1}^{N-k} \psi^i_s \psi^{\bar j}_{\bar s} \psi^j_t \psi^{\bar i}_{\bar t} \right).
\end{align}
In the second term, $( \omega^{ii}\omega^{jj})^{N-k-1}$ contains $N-k-1$ $\psi^i_s \psi^{\bar i}_{\bar s}$, and $\psi^j_t \psi^{\bar i}_{\bar t}$. Therefore, it must be $s=t$ based on the conditions of the fermionic variables.
\begin{align}
&P_1= (kN-k^2-2)!k(N-k)(N-k-1) \notag \\
&-\sum_{i\neq j} \sum_{s=1}^{N-k} \frac{(kN-k^2-2)!}{((N-k-1)!)^2} \int D\psi \left( \prod_{n\neq i,j} \prod_{l=1}^{N-k} \psi^n_l \psi^{\bar n}_{\bar l} \right) \left( \omega^{ii}\omega^{jj}\right)^{N-k-1} \left(  \psi^i_s \psi^{\bar i}_{\bar s} \psi^j_s \psi^{\bar j}_{\bar s}   \right)\\
&=(kN-k^2-2)!k(N-k)(N-k-1) \notag \\
&-\sum_{i\neq j} \sum_{s=1}^{N-k} (kN-k^2-2)! \int D\psi \left( \prod_{n\neq i,j} \prod_{l=1}^{N-k} \psi^n_l \psi^{\bar n}_{\bar l} \right) \left( \prod_{\substack{q=1 \\q\neq s}}^{N-k}\psi^i_q \psi^{\bar i}_{\bar q} \psi^j_q \psi^{\bar j}_{\bar q}  \right) \left(  \psi^i_s \psi^{\bar i}_{\bar s} \psi^j_s \psi^{\bar j}_{\bar s}   \right)\\
&=(kN-k^2-2)!k(N-k)(N-k-1)-\sum_{i\neq j} \sum_{s=1}^{N-k} (kN-k^2-2)!  \\
&=(kN-k^2-2)!\{k(N-k)(N-k-1)-(N-k)k(k-1)\}=(kN-k^2-2)!k(N-k)(N-2k).
\end{align}
Therefore, we obtain $P_1$. We compute $P_2$.
\begin{align}
&P_2=\int D\psi \left( \sum_{n=1}^k \omega^{nn} \right)^{kN-k^2-4} \left(\sum_{i=1}^k (\omega^{ii})^2+\sum_{i\neq j} \omega^{ij}\omega^{ji} \right)^2 \\
&=\int D\psi \left( \sum_{n=1}^k \omega^{nn} \right)^{kN-k^2-4} \left[\sum_{i,j} (\omega^{ii}  \omega^{jj})^2+2\sum_{m=1}^k\sum_{i\neq j} (\omega^{mm})^2\omega^{ij}\omega^{ji} +\sum_{a\neq b}\sum_{i\neq j} \omega^{ab}\omega^{ba}\omega^{ij}\omega^{ji} \right]. 
\end{align}
We define $Q_1$, $Q_2$ and $Q_3$ as follows.
\begin{align}
Q_1&:=\sum_{i,j}\int D\psi \left( \sum_{n=1}^k \omega^{nn} \right)^{kN-k^2-4}  (\omega^{ii}  \omega^{jj})^2,
Q_2:=2\sum_{m=1}^k\sum_{i\neq j} \int D\psi \left( \sum_{n=1}^k \omega^{nn} \right)^{kN-k^2-4}(\omega^{mm})^2\omega^{ij}\omega^{ji}, \\
Q_3&:=\sum_{a\neq b}\sum_{i\neq j} \int D\psi \left( \sum_{n=1}^k \omega^{nn} \right)^{kN-k^2-4}\omega^{ab}\omega^{ba}\omega^{ij}\omega^{ji}.
\end{align}
Thereafter, $P_2=Q_1+Q_2+Q_3$. We consider $Q_1$:  
\begin{align}
Q_1=\sum_{i=1}^k\int D\psi \left( \sum_{n=1}^k \omega^{nn} \right)^{kN-k^2-4}  (\omega^{ii} )^4
+\sum_{i\neq j}\int D\psi \left( \sum_{n=1}^k \omega^{nn} \right)^{kN-k^2-4}  (\omega^{ii}  \omega^{jj})^2.
\end{align}
 We can compute the above equation in the same manner as $P_1$. Consequently,
\begin{align}
Q_1=(kN-k^2-4)!k(N-k)\{(N-k-1)(N-k-2)(N-k-3)+(k-1)(N-k)(N-k-1)^2\}.
\end{align}
Subsequently, we calculate $Q_2$.
\begin{align}
Q_2&=2\sum_{i\neq j} \int D\psi \left( \sum_{n=1}^k \omega^{nn} \right)^{kN-k^2-4}(\omega^{ii})^2\omega^{ij}\omega^{ji} 
+2\sum_{i\neq j} \int D\psi \left( \sum_{n=1}^k \omega^{nn} \right)^{kN-k^2-4}(\omega^{jj})^2\omega^{ij}\omega^{ji} \notag \\
&+2\sum_{i\neq j}\sum_{m\neq i,j} \int D\psi \left( \sum_{n=1}^k \omega^{nn} \right)^{kN-k^2-4}(\omega^{mm})^2\omega^{ij}\omega^{ji}.
\end{align}
From $\omega^{ij}\omega^{ji} =\omega^{ji} \omega^{ij}$, if we replace $i$ with $j$ and $j$ with $i$ in the second term, it is the same as in the first term.  
\begin{align}
Q_2&=4\sum_{i\neq j} \int D\psi \left( \sum_{n=1}^k \omega^{nn} \right)^{kN-k^2-4}(\omega^{ii})^2\omega^{ij}\omega^{ji} +2\sum_{i\neq j}\sum_{m\neq i,j} \int D\psi \left( \sum_{n=1}^k \omega^{nn} \right)^{kN-k^2-4}(\omega^{mm})^2\omega^{ij}\omega^{ji} \\
&=4\sum_{i\neq j} \sum_{p_n}\frac{(kN-k^2-4)!}{\prod_{q=1}^k p_q!} \int D\psi \left( \prod_{n=1}^k (\omega^{nn})^{p_n} \right)(\omega^{ii})^2\omega^{ij}\omega^{ji} \notag \\
&+2\sum_{i\neq j}\sum_{m\neq i,j} \sum_{p_n}\frac{(kN-k^2-4)!}{\prod_{q=1}^k p_q!} \int D\psi \left( \prod_{n=1}^k (\omega^{nn})^{p_n} \right)(\omega^{mm})^2\omega^{ij}\omega^{ji}.
\end{align}
Thereafter, $\sum_{p_n}$ is the sum of all combinations from $0$ to $kN-k^2-4$ indices $p_1$ through $p_k$ such that the sum of  $p_n (n=1,\cdots,k)$ is $kN-k^2-4$. From the condition of fermionic integration and the condition of  fermionic variables $\psi^i_s \psi^i_s=0$, in the first term, $p_n=N-k (n\neq i,j)$ and $p_i=N-k-3$, $p_j=N-k-1$. In the second term, $p_n=N-k (n\neq i,j,m)$ and $p_i=p_j=N-k-1$, $p_m=N-k-2$. Therefore, 
\begin{align}
Q_2
&=4\sum_{i\neq j} \frac{(kN-k^2-4)!}{(N-k-3)!(N-k-1)!} \int D\psi \left( \prod_{n\neq i,j} \prod_{l=1}^{N-k}  \psi^n_l \psi^{\bar n}_{\bar l} \right)(\omega^{ii}\omega^{jj})^{N-k-1}\omega^{ij}\omega^{ji} \notag \\
&+2\sum_{i\neq j}\sum_{m\neq i,j} \frac{(kN-k^2-4)!(N-k)!}{(N-k-2)!((N-k-1)!)^2} \int D\psi \left( \prod_{n\neq i, j} \prod_{l=1}^{N-k}  \psi^n_l \psi^{\bar n}_{\bar l} \right)(\omega^{ii}\omega^{jj})^{N-k-1}\omega^{ij}\omega^{ji}.
\end{align}
Here,  we can  calculate the fermion integral in the same manner as $P_1$. We obtain
\begin{align}
\int D\psi &\left( \prod_{n\neq i, j} \prod_{l=1}^{N-k}  \psi^n_l \psi^{\bar n}_{\bar l} \right)(\omega^{ii}\omega^{jj})^{N-k-1}\omega^{ij}\omega^{ji}=-(N-k)((N-k-1)!)^2. \\
Q_2
&=-4\sum_{i\neq j} \frac{(kN-k^2-4)!(N-k)!}{(N-k-3)!} -2\sum_{i\neq j} \sum_{m\neq i,j}\frac{(kN-k^2-4)!(N-k)!(N-k)}{(N-k-2)!} \\
&=-4\frac{(kN-k^2-4)!(N-k)!}{(N-k-3)!}k(k-1) -2\frac{(kN-k^2-4)!(N-k)!(N-k)}{(N-k-2)!} k(k-1)(k-2) \notag \\
&=(kN-k^2-4)!k(N-k)(k-1)[-4(N-k-1)(N-k-2)-2(N-k)(N-k-1)(k-2)].
\end{align}
Finally, we compute $Q_3$.
\begin{align}
&Q_3=\sum_{a\neq b}\sum_{i\neq j} \int D\psi \left( \sum_{n=1}^k \omega^{nn} \right)^{kN-k^2-4}\omega^{ab}\omega^{ba}\omega^{ij}\omega^{ji}.
\end{align}
The sum $\sum_{a\neq b}\sum_{i\neq j}$ can be divided into the following seven cases.
\begin{itembox}{Sum patterns of $(i,j)$ and $(a,b)$}

(1) $i=a$, $j=b$.  (2) $i=b$, $j=a$. (3) $i=a$, $j\neq b$. (4) $i=b$, $j\neq a$. (5) $i\neq a$, $j=b$. (6) $i\neq b$, $j=a$. (7) $i\neq a,b$ $j\neq a,b$.
\end{itembox}
From the symmetry of $a$, $b$ and $i$, $j$, (1) and (2) have the same form: Similarly, (3), (4), (5), and (6) have the same form: Therefore, 
\begin{align}
Q_3&=2\sum_{i\neq j} \int D\psi \left( \sum_{n=1}^k \omega^{nn} \right)^{kN-k^2-4}(\omega^{ij}\omega^{ji})^2
+4\sum_{i\neq j} \sum_{b\neq i,j}  \int D\psi \left( \sum_{n=1}^k \omega^{nn} \right)^{kN-k^2-4}\omega^{ib}\omega^{bi}\omega^{ij}\omega^{ji} \notag \\
&+\sum_{(i,j,a,b)}^\prime  \int D\psi \left( \sum_{n=1}^k \omega^{nn} \right)^{kN-k^2-4}\omega^{ab}\omega^{ba}\omega^{ij}\omega^{ji}.
\end{align}
Here, $\sum_{(i,j,a,b)}^\prime$ implies  that $i, j, a$, and $b$ are different from each other in the summation.
\begin{align}
Q_3&=2\sum_{i\neq j} \frac{(kN-k^2-4)!}{((N-k)!)^{k-2}((N-k-2)!)^2}\int D\psi \left( \prod_{n\neq i,j} (\omega^{nn})^{N-k} \right)\left( \omega^{ii} \omega^{jj} \right)^{N-k-2} (\omega^{ij}\omega^{ji})^2 \notag \\
&+4\sum_{i\neq j} \sum_{b\neq i,j} \frac{(kN-k^2-4)!}{((N-k)!)^{k-3}(N-k-2)!((N-k-1)!)^2}  \notag \\
&\times \int D\psi \left( \prod_{n\neq i,j,b} (\omega^{nn})^{N-k} \right) \left(\omega^{ii} \right)^{N-k-2} \left( \omega^{bb} \omega^{jj} \right)^{N-k-1} \omega^{ib}\omega^{bi}\omega^{ij}\omega^{ji} \notag \\
&+\sum_{(i,j,a,b)}^\prime \frac{(kN-k^2-4)!}{((N-k)!)^{k-4} ((N-k-1)!)^{4} } \int D\psi \left( \prod_{n\neq a,b,i,j} (\omega^{nn})^{N-k} \right) \left(\omega^{aa} \omega^{bb}\omega^{ii} \omega^{jj} \right)^{N-k-1} \notag \\
&\times \omega^{ab}\omega^{ba}\omega^{ij}\omega^{ji} \\
&=2\sum_{i\neq j} \frac{(kN-k^2-4)!}{((N-k-2)!)^2}\int D\psi \left( \prod_{n\neq i,j} \prod_{l=1}^{N-k}\psi^n_l \psi^{\bar n}_{\bar l} \right)\left( \omega^{ii} \omega^{jj} \right)^{N-k-2} (\sum_{s_1,s_2,t_1,t_2}\psi^i_{s_1} \psi^{\bar j}_{\bar s_1} \psi^j_{t_1} \psi^{\bar i}_{\bar t_1} \psi^i_{s_2} \psi^{\bar j}_{\bar s_2} \psi^j_{t_2} \psi^{\bar i}_{\bar t_2} ) \notag \\
&+4\sum_{i\neq j} \sum_{b\neq i,j} \frac{(kN-k^2-4)!}{(N-k-2)!((N-k-1)!)^2}  \notag \\
&\times \int D\psi \left( \prod_{n\neq i,j} \prod_{l=1}^{N-k}\psi^n_l \psi^{\bar n}_{\bar l} \right) \left(\omega^{ii} \right)^{N-k-2} \left( \omega^{bb} \omega^{jj} \right)^{N-k-1} (\sum_{s_1,s_2,t_1,t_2}\psi^i_{s_1} \psi^{\bar b}_{\bar s_1} \psi^b_{t_1} \psi^{\bar i}_{\bar t_1} \psi^i_{s_2} \psi^{\bar j}_{\bar s_2} \psi^j_{t_2} \psi^{\bar i}_{\bar t_2} ) \notag \\
&+\sum_{(i,j,a,b)}^\prime \frac{(kN-k^2-4)!}{ ((N-k-1)!)^{4} } \int D\psi \left( \prod_{n\neq a,b,i,j} \prod_{l=1}^{N-k}\psi^n_l \psi^{\bar n}_{\bar l} \right) \left(\omega^{aa} \omega^{bb}\omega^{ii} \omega^{jj} \right)^{N-k-1} \notag \\&(\sum_{s_1,s_2,t_1,t_2}\psi^a_{s_1} \psi^{\bar b}_{\bar s_1} \psi^b_{t_1} \psi^{\bar a}_{\bar t_1} \psi^i_{s_2} \psi^{\bar j}_{\bar s_2} \psi^j_{t_2} \psi^{\bar i}_{\bar t_2} ).
\end{align}
We consider sums of $s_1,s_2,t_1$ and $t_2$. In the first term, the sum can be divided into two ways, ($s_1=t_1$, $s_2=t_2$, $s_1 \neq s_2$) and ($s_1=t_2$, $s_2=t_1$,$s_1 \neq s_2$). In the second term, it must be ($s_1=t_1$, $s_2=t_2$, $s_1\neq s_2$). In the third term, it must be ($s_1=t_1$, $s_2=t_2$). Because the first term is symmetric for $s_1$ and $s_2$ and $t_1$ and $t_2$,
\begin{align}
Q_3
&=4\sum_{i\neq j}\sum_{s_1\neq s_2} \frac{(kN-k^2-4)!}{((N-k-2)!)^2}\int D\psi \left( \prod_{n\neq i,j} \prod_{l=1}^{N-k}\psi^n_l \psi^{\bar n}_{\bar l} \right)\left( \omega^{ii} \omega^{jj} \right)^{N-k-2} (\psi^i_{s_1} \psi^{\bar i}_{\bar s_1}  \psi^j_{s_1} \psi^{\bar j}_{\bar s_1} \psi^i_{s_2}  \psi^{\bar i}_{\bar s_2} \psi^j_{s_2}  \psi^{\bar j}_{\bar s_2} ) \notag \\
&+4\sum_{i\neq j} \sum_{b\neq i,j} \sum_{s_1\neq s_2}\frac{(kN-k^2-4)!}{(N-k-2)!((N-k-1)!)^2}  \notag \\
&\times \int D\psi \left( \prod_{n\neq i,j} \prod_{l=1}^{N-k}\psi^n_l \psi^{\bar n}_{\bar l} \right) \left(\omega^{ii} \right)^{N-k-2} \left( \omega^{bb} \omega^{jj} \right)^{N-k-1} (\psi^i_{s_1} \psi^{\bar i}_{\bar s_1} \psi^b_{s_1} \psi^{\bar b}_{\bar s_1}   \psi^i_{s_2} \psi^{\bar i}_{\bar s_2}  \psi^j_{s_2} \psi^{\bar j}_{\bar s_2} ) \notag \\
&+\sum_{(i,j,a,b)}^\prime \sum_{s_1,s_2} \frac{(kN-k^2-4)!}{ ((N-k-1)!)^{4} } \int D\psi \left( \prod_{n\neq a,b,i,j} \prod_{l=1}^{N-k}\psi^n_l \psi^{\bar n}_{\bar l} \right) \left(\omega^{aa} \omega^{bb}\omega^{ii} \omega^{jj} \right)^{N-k-1} \notag \\&(\psi^a_{s_1} \psi^{\bar a}_{\bar s_1} \psi^b_{s_1} \psi^{\bar b}_{\bar s_1}  \psi^i_{s_2} \psi^{\bar i}_{\bar s_2}  \psi^j_{s_2} \psi^{\bar j}_{\bar s_2} )  \\
&=4\sum_{i\neq j}\sum_{s_1\neq s_2} (kN-k^2-4)!+4\sum_{i\neq j} \sum_{b\neq i,j} \sum_{s_1\neq s_2}(kN-k^2-4)! +\sum_{(i,j,a,b)}^\prime \sum_{s_1,s_2} (kN-k^2-4)! \\
&=(kN-k^2-4)![4k(k-1)(N-k)(N-k-1)+4k(k-1)(k-2)(N-k)(N-k-1) \notag \\
&+k(k-1)(k-2)(k-3)(N-k)^2] \\
&=(kN-k^2-4)!k(N-k)[4(k-1)^2(N-k-1) +(k-1)(k-2)(k-3)(N-k)].
\end{align}
Therefore, we obtain (\ref{prop1-2}) from $P_3=Q_1+Q_2+Q_3$ and these results. We complete the proof of proposition 1.
\end{proof}

\section*{Acknowledgement}
We would like thank Professor Masao Jinzenji for the useful discussions. This paper would not have been completed without his help.

\bibliographystyle{plain}
\bibliography{KK-revised5}

\begin{thebibliography}{10}

\bibitem{AW}
A.Weber.
\newblock Equivariant {C}hern classes and localization theorem.
\newblock {\em J. Singul.}, 5:153--176, 2012.
\newblock {DOI}:10.5427/jsing.2012.5k, {URL}:
  https://doi.org/10.5427/jsing.2012.5k. arXiv:1110.5515.

\bibitem{DTH}
D.T.Hiep.
\newblock Identities involving (doubly) symmetric polynomials and integrals
  over {G}rassmannians, 2016.
\newblock arXiv, {DOI}: 10.48550/ARXIV.1607.04850 {URL}:
  https://arxiv.org/abs/1607.04850.

\bibitem{EW3}
E.Witten.
\newblock The {V}erlinde algebra and the cohomology of the {G}rassmannian.
\newblock In {\em Geometry, topology, \& physics}, Conf. Proc. Lecture Notes
  Geom. Topology, IV, pages 357--422. Int. Press, Cambridge, MA, 1995.
\newblock arXiv:hep-th/9312104.

\bibitem{MZ}
M.Zielenkiewicz.
\newblock Integration over homogeneous spaces for classical {L}ie groups using
  iterated residues at infinity.
\newblock {\em centr.eur.j.math.}, 12:574--583, 2014.
\newblock {DOI}:10.2478/s11533-013-0372-z, {URL}:
  https://doi.org/10.2478/s11533-013-0372-z. arXiv:1212.6623.

\bibitem{NC}
N.Chair.
\newblock Intersection numbers on {G}rassmannians, and on the space of
  holomorphic maps from ${C P}^1$ into ${G}_r( {C}^n)$.
\newblock {\em J. Geom. Phys.}, 38(2):170--182, 2001.
\newblock {DOI}:10.1016/S0393-0440(00)00059-0,
  {URL}:https://doi.org/10.1016/S0393-0440(00)00059-0. arXiv:hep-th/9808170.

\bibitem{GH}
P.Griffiths and J.Harris.
\newblock {\em Principles of algebraic geometry}.
\newblock Pure and Applied Mathematics. Wiley--Interscience, New York, 1978.

\bibitem{BT}
R.Bott and L.W.Tu.
\newblock {\em Differential Forms in Algebraic Topology}.
\newblock Graduate Texts in Mathematics. Springer New York, NY, 1982.

\bibitem{Shoichiro}
S.Imanishi.
\newblock Computation of {E}uler number for the {G}rassmann manifold {G}(2,{N})
  via {M}athai-{Q}uillen formalism (in {J}apanese, printed in {J}apan).
\newblock Master's thesis, Hokkaido University, 2019.

\bibitem{IJK}
S.Imanishi, M.Jinzenji, and K.Kuwata.
\newblock Evaluation of {E}uler number of complex {G}rassmann manifold
  {G}(k,{N}) via {M}athai-{Q}uillen formalism.
\newblock {\em Journal of Geometry and Physics}, 180:104623, 2022.
\newblock {DOI}:https://doi.org/10.1016/j.geomphys.2022.104623.

\bibitem{IN}
T.Ikeda and H.Naruse.
\newblock Modern {S}chubert calculus, from the special polynomial theory's
  point of view (in {J}apanese, printed in {J}apan).
\newblock {\em Sugaku}, 63(3):313--337, 2011.
\newblock {DOI}:10.11429/sugaku.0633313,
  {URL}:https://doi.org/10.11429/sugaku.0633313.

\bibitem{WF}
W.Fulton.
\newblock {\em Intersection {T}heory}.
\newblock Springer-Verlag Berlin Hidelberg, 1998.

\end{thebibliography}

\end{document}